\documentclass[12pt,reqno,a4wide]{amsart}

\usepackage{tikz} 
\usepackage{calrsfs}

\allowdisplaybreaks

\newcommand{\ga}[2]{\genfrac{[}{]}{0pt}{}{#1}{#2}}

\newcommand{\gaa}[2]{\genfrac{[}{]}{0pt}{}{#1}{#2}_{q^2}}

\title{A new recursion for Bressoud's polynomials}

\author{Helmut Prodinger}
\address{Helmut Prodinger\\
	Department of Mathematical Sciences\\
	Stellenbosch University\\
	7602 Stellenbosch\\
	South Africa}
\email{hproding@sun.ac.za}

\begin{document}
\begin{abstract}
	A new recursion in only one variable allows very simple verifications of Bressoud's polynomial identities, which lead to the Rogers-Ramanujan identities. This approach might be compared with an earlier approach due to Chapman. Applying the
	$q$-Chu-Vandermonde convolution, as suggested by Cigler, makes the computations particularly simple and elementary.
	The same treatment is also applied to the Santos polynomials and perhaps more polynomials from a list of Rogers-Ramanujan like polynomials~\cite{Sills}.
	\end{abstract}

\maketitle
	
\section{Bressoud's polynomial identities}

Let
\begin{equation*}
A_n=\sum_{k=0}^nq^{k^2}\ga{n}{k},\qquad B_n=\sum_{j\in\mathbb{Z}}(-1)^jq^{\frac{j(5j-1)}{2}}\ga{2n}{n-2j},
\end{equation*}	
and
\begin{equation*}
C_n=\sum_{k=0}^nq^{k^2+k}\ga{n}{k},\qquad D_n=\sum_{j\in\mathbb{Z}}(-1)^jq^{\frac{j(5j-3)}{2}}\ga{2n+1}{n+1-2j},
\end{equation*}	
where $\ga nk$ is a $q$-binomial coefficients~\cite{Andrews-Eriksson}, defined by
\begin{equation*}
\ga nk:=\frac{(q;q)_{n}}{(q;q)_{k}(q;q)_{n-k}}\quad\text{with}\quad (x;q)_m:=(1-x)(1-xq)\dots(1-xq^{m-1}).
\end{equation*}
	
Bressoud~\cite{Bressoud} proved that $A_n=B_n$ and $C_n=D_n$ and that taking the limit $n\to\infty$	leads to the 
celebrated Rogers-Ramanujan identities. Since it is well documented in the literature how to take this limit we will not
repeat this here and concentrate on the polynomial identities. The Bressoud polynomials $A_n$ and $C_n$ are not the only finitizations of the celebrated Rogers-Ramanujan identities, but arguably the simplest and prettiest. More information about finite
versions of Rogers-Ramanujan type identities can be found in the encyclopedic paper~\cite{Sills}.

Chapman~\cite{Chapman} found a simple and elementary approach, and it was  used in \cite{Andrews-Eriksson} almost without change. A different simple proof was provided by Cigler~\cite{Cigler} a few years later.

Chapman's method consists in showing that both sides of the identity satisfy the same recursion. This recursion is, however, in two variables. Also, auxiliary sequences needed to be introduced. 

Here, we use a different (although related) recursion that depends only on one variable, and requires no auxiliary sequences.

It is easy to check that the first two values of the sequences also coincide, so that the sequences themselves coincide.

There are other approaches to deal with Bressoud's polynomials and extensions, like \cite{Warnaar}. Here, we try to make everything as simple and elementary as possible.

In a final section, we apply the same machinery to the so-called Santos polynomials \cite{AndrewsSantos}. They belong to another pair of Rogers-Ramanujan type identities, and there is hope that even more such polynomials can be treated along the lines of this note, since Zeilberger's algorithm is helpful to establish the relevant recursions.

\section{The first identity}

It is a routine computation to verify that
\begin{equation}\label{last}
\ga{n}{k}-(1+q-q^n)\ga{n-1}{k}-q^{2n-2k}\ga{n-1}{k-1}+q(1-q^{n-1})\ga{n-2}{k}=0.
\end{equation} 
Multiplying this by $q^{k^2}$ and summing over nonnegative integers $k$ leads to the recursion
\begin{equation*}
A_n-(1+q-q^n+q^{2n-1})A_{n-1}+q(1-q^{n-1})A_{n-2}=0.
\end{equation*}
It is more of a challenge to show the recursion	
 \begin{equation}\label{rec-B}
B_n-(1+q-q^n+q^{2n-1})B_{n-1}+q(1-q^{n-1})B_{n-2}=0,
\end{equation}
which we will do now.

Cigler~\cite{Cigler} used the $q$-Chu-Vandermonde convolution to  expand
 \begin{equation*}
\ga{2n}{n-2j}= \sum_kq^{(k-j)(k+j)}\ga{n}{k-j}\ga{n}{k+j}
 \end{equation*}
and consequently
 	\begin{equation*}
 	B_n=\sum_kq^{k^2}f(n,k)
 	\end{equation*}
with
 \begin{align*}
f(n,k):=\sum_{j}(-1)^jq^{\frac{j(3j-1)}{2}}\ga{n}{k-j}\ga{n}{k+j}.
 \end{align*}
This identity has (of course) a long history. As one referee has pointed out, this identity is at the core of the analysis of Durfee squares~\cite{Andrews76}, but can also be derived from classical $q$-hypergeometric identities ($q$-Dougall, $q$-Dixon, \dots), as already pointed out by Cigler himself, with a reference to Warnaar~\cite{Warnaar}.  It is, however, quite striking, that this polynomial identity can be used in elementary proofs related to Bressoud's polynomials.

 One can actually compute $f(n,k)=\ga{n}{k}$, but this will play no role in our proof. We only need simple recursions for $f(n,k)$ that
 appear already in \cite{Cigler} but are repeated here for completeness. For that, only the recursions for the $q$-binomial coefficients are needed. (These recursions are the two standard recursions that prove $f(n,k)=\ga{n}{k}$ anyway.)
 
 We start with
 \begin{align*}
\sum_j&(-1)^jq^{\frac{j(3j-1)}{2}}\ga{n-1}{k-j}\ga{n}{k+j}\\
&=\sum_j(-1)^jq^{\frac{j(3j-1)}{2}}\ga{n-1}{k-j}\ga{n-1}{k+j}+q^{n-k}\sum_j(-1)^jq^{\frac{3j(j-1)}{2}}\ga{n-1}{k-j}\ga{n-1}{k+j-1}\\
&=f(n-1,k),
 \end{align*}
 since the last sum, on the substitution $j\to-j+1$, turns into its own negative.
 Similarly,
 \begin{align*}
 \sum_j&(-1)^jq^{\frac{j(3j+1)}{2}}\ga{n-1}{k-j}\ga{n}{k+j+1}\\
 &=q^{k+1}\sum_j(-1)^jq^{\frac{3j(j+1)}{2}}\ga{n-1}{k-j}\ga{n-1}{k+j+1}+\sum_j(-1)^jq^{\frac{j(3j+1)}{2}}\ga{n-1}{k-j}\ga{n-1}{k+j}\\
 &=f(n-1,k),
 \end{align*}
 by the same reasoning. Therefore
 \begin{align*}
f(n,k)&=\sum_{j}(-1)^jq^{\frac{j(3j-1)}{2}}\ga{n}{k+j}\bigg(\ga{n-1}{k-j}+q^{n-k+j}\ga{n-1}{k-j-1}\bigg)\\
&=f(n-1,k)+q^{n-k}\sum_{j}(-1)^jq^{\frac{j(3j+1)}{2}}\ga{n}{k+j}\ga{n-1}{k-j-1}\\
&=f(n-1,k)+q^{n-k}f(n-1,k-1).
\end{align*}
A very similar computation leads to
\begin{align*}
f(n,k)&=\sum_{j}(-1)^jq^{\frac{j(3j-1)}{2}}\ga{n}{k-j}\bigg(q^{k+j}\ga{n-1}{k+j}+\ga{n-1}{k+j-1}\bigg)\\
&=q^k\sum_{j}(-1)^jq^{\frac{j(3j+1)}{2}}\ga{n}{k-j}\ga{n-1}{k+j}+\sum_{j}(-1)^jq^{\frac{j(3j-1)}{2}}\ga{n}{k-j}\ga{n-1}{k+j-1}\\
&=q^k\sum_{j}(-1)^jq^{\frac{j(3j-1)}{2}}\ga{n}{k+j}\ga{n-1}{k-j}+\sum_{j}(-1)^jq^{\frac{j(3j+1)}{2}}\ga{n}{k+j}\ga{n-1}{k-j-1}\\
&=q^kf(n-1,k)+f(n-1,k-1).
\end{align*}
Now we can prove that the sequence $B_n$  satisfies the recursion (\ref{rec-B}), which means that
\begin{align*}
\sum_kq^{k^2}f(n,k)&-(1+q-q^n)\sum_kq^{k^2}f(n-1,k)-q^{2n}\sum_kq^{k^2-2k}f(n-1,k-1)\\&+q(1-q^{n-1})\sum_kq^{k^2}f(n-2,k)=0.
\end{align*} 
 We claim that even
\begin{align*}
f(n,k)-(1+q-q^n)f(n-1,k)-q^{2n-2k}f(n-1,k-1)+q(1-q^{n-1})f(n-2,k)=0.
\end{align*} 
Using the recursion \ref{{last}}, this is equivalent to
\begin{align*}
&q^{n-k}f(n-1,k-1)-(q-q^n)f(n-1,k)-q^{2n-2k}f(n-1,k-1)+q(1-q^{n-1})f(n-2,k)\\
&=(q^{n-k}-q^{2n-2k})f(n-1,k-1)-(q-q^n)(f(n-1,k)-f(n-2,k))\\
&=(q^{n-k}-q^{2n-2k})f(n-1,k-1)-(q-q^n)q^{n-k-1}f(n-2,k-1)=0.
\end{align*} 
This may be further reduced to
\begin{align*}
(1-&q^{n-k})f(n-1,k-1)-(1-q^{n-1})f(n-2,k-1)\\
&=f(n-1,k-1)-f(n-2,k-1)-q^{n-k}(f(n-1,k-1)-q^{k-1}f(n-2,k-1))\\
&=q^{n-k}f(n-2,k-2)-q^{n-k}f(n-2,k-2)=0,
\end{align*} 
which is now obvious.

\textbf{Remark.} As one referee has pointed out, the two coupled recursions appearing in \cite{Chapman} can be transformed into the recursion (\ref{rec-B}). This goes as follows, starting from
\begin{align*}
B_n&=B_{n-1}+q^nD_{n-1},\\
D_n-q^nB_n&=(1-q^n)D_{n-1}.
\end{align*}
Eliminating $D_{n-1}$ from the second equation and replacing it in the first equation leads to
\begin{align*}
(1-q^n)B_n&=(1-q^n)B_{n-1}+q^nD_n-q^{2n}B_n\\
&=(1-q^n)B_{n-1}+\frac1q(B_{n+1}-B_n)-q^{2n}B_n.
\end{align*}
Rearranging this leads to
\begin{equation*}
B_{n+1}=(q-q^{n+1}+1+q^{2n+1})B_n-q(1-q^n)B_{n-1}.
\end{equation*}
Replacing $n$ by $n-1$ leads to the recursion (\ref{rec-B}). 

The approach in the present paper is to find the second order recursion (in only one variable) \emph{directly}, which will we used in the sequel for the second identity, as well as for the Santos-polynomials. Fortunately, the $q$-Zeilberger algorithm helps to find it if one does not see it otherwise. In the instance of Santos-polynomials, the recursions (\ref{santos-rec}), (\ref{santosT-rec}) are readily found with a computer, but an elimination using auxiliary sequences would be a more elaborate process.

\section{The second identity}

From
\begin{equation*}
\ga{n}{k}-(1+q-q^n)\ga{n-1}{k}-q^{2n-2k}\ga{n-1}{k-1}+q(1-q^{n-1})\ga{n-2}{k}=0,
\end{equation*} 
multiplying this by $q^{k^2+k}$ and summing over all nonnegative integers $k$ we are led to the recursion
\begin{equation*}
C_n-(1+q-q^n+q^{2n})C_{n-1}+q(1-q^{n-1})C_{n-2}=0.
\end{equation*}
Now we will deduce the recursion
\begin{equation*}
D_n-(1+q-q^n+q^{2n})D_{n-1}+q(1-q^{n-1})D_{n-2}=0
\end{equation*}
as well.

The $q$-Chu-Vandermonde formula leads to
\begin{equation*}
\ga{2n+1}{n-2j}=\sum_k q^{k^2+k-j^2+j}\ga{n+1}{k+1-j}\ga{n}{k+j},
\end{equation*}
and therefore
\begin{align*}
D_n&=\sum_{j}(-1)^jq^{\frac{j(5j-3)}{2}}\sum_k q^{k^2+k-j^2+j}\ga{n+1}{k+1-j}\ga{n}{k+j}\\
&=\sum_k q^{k^2+k}\sum_{j}(-1)^jq^{\frac{j(3j-1)}{2}}\ga{n+1}{k+1-j}\ga{n}{k+j}\\
&=\sum_k q^{k^2+k}f(n,k).
\end{align*}
From
\begin{align*}
	f(n,k)-(1+q-q^n)f(n-1,k)-q^{2n-2k}f(n-1,k-1)+q(1-q^{n-1})f(n-2,k)=0
\end{align*} 
we get, upon multiplication with $q^{k^2+k}$ and summing over all nonnegative integers $k$
\begin{align*}
D_n-(1+q-q^n)D_{n-1}-q^{2n}D_{n-1}+q(1-q^{n-1})D_{n-2}=0,
\end{align*} 
as claimed.

\section{Santos polynomials}

The Santos polynomials are defined as
\begin{equation*}
S_n:=\sum_{0\le 2k\le n}q^{2k^2}\ga{n}{2k}.
	\end{equation*}
	The are used to prove identities A.39 and A.38 from the list \cite{Sills}.
	
We start from
\begin{equation}\label{santos-rec}
\ga{n+2}{2k}-(1+q)\ga{n+1}{2k}+q\ga{n}{2k}-q^{2n+4-4k}\ga{n}{2k-2}=0.
\end{equation}
Multiplying this by $q^{2k^2}$ and summing, we find the recursion
\begin{equation*}
S_{n+2}-(1+q)S_{n+1}+(q-q^{2n+2})S_n=0.
\end{equation*}
(Originally, it was found using Zeilberger's $q$-EKHAD algorithm.)

The alternative form for the Santos polynomials, as to be shown, is
\begin{equation*}
\overline{S}_n=\sum_{j}q^{4j^2-j}\gaa{n}{\lfloor\tfrac{n+1}{2} \rfloor-2j}.
\end{equation*}
The $q$-Chu-Vandermonde formula leads to
\begin{equation*}
\gaa{n}{\lfloor\tfrac{n+1}{2} \rfloor-2j}=\sum_k \gaa{\lceil\tfrac{n}{2} \rceil}{k+j}\gaa{\lfloor\tfrac{n}{2} \rfloor}{k-j}
q^{2k^2-2j^2}.
\end{equation*}
Therefore
\begin{align*}
\overline{S}_n&=\sum_kq^{2k^2}g(n,j)
\end{align*}
with
\begin{equation*}
g(n,j)=\sum_{j}q^{2j^2-j} \gaa{\lceil\tfrac{n}{2} \rceil}{k+j}\gaa{\lfloor\tfrac{n}{2} \rfloor}{k-j}.
\end{equation*}
Using only the recursions for the $q$-binomial coefficients, a tedious computation leads to
\begin{equation*}
	g(n+2,j)-(1+q)g(n+1,j)+q\,g(n,j)-q^{2n+2-2j}g(n,j-1)=0.
\end{equation*}
Multiplying this by $q^{2k^2}$ and summing over all nonnegative integers $k$, we find the recursion
\begin{equation*}
	\overline{S}_{n+2}-(1+q)\overline{S}_{n+1}+(q-q^{2n+2})\overline{S}_n=0.
\end{equation*}
Since the recursion for $g(n,j)$ defines, together with some initial conditions, the sequence uniquely,
this also shows that $g(n,j)=\ga n{2j}$.

There is a second family of Santos polynomials, defined by
\begin{equation*}
	T_n:=\sum_{0\le 2k+1\le n}q^{2k^2+2k}\ga{n}{2k+1}.
\end{equation*}

We start from
\begin{equation}\label{santosT-rec}
	\ga{n+2}{2k+1}-(1+q)\ga{n+1}{2k+1}+q\ga{n}{2k+1}-q^{2n+2-4k}\ga{n}{2k-1}=0.
\end{equation}
Multiplying this by $q^{2k^2+2k}$ and summing leads to
\begin{equation*}
	T_{n+2}-(1+q)T_{n+1}+(q-q^{2n+2})T_n=0.
\end{equation*}

The alternative form for the second family of Santos polynomials, as to be shown, is
\begin{equation*}
	\overline{T}_n=\sum_{j}q^{4j^2-3j}\gaa{n}{\lfloor\tfrac{n+2}{2} \rfloor-2j}.
\end{equation*}
The $q$-Chu-Vandermonde formula leads to
\begin{equation*}
	\gaa{n}{\lfloor\tfrac{n+2}{2} \rfloor-2j}=\sum_k \gaa{\lfloor\tfrac{n}{2} \rfloor}{k+j}\gaa{\lceil\tfrac{n}{2} \rceil}{k+1-j}
	q^{2k^2+2k-2j^2+2j}
\end{equation*}
and therefore
\begin{equation*}
	\overline{T}_n=\sum_kq^{2k^2+2k}h(n,k),
\end{equation*}
with
\begin{equation*}
h(n,k):=\sum_{j}q^{2j^2-j} \gaa{\lfloor\tfrac{n}{2} \rfloor}{k+j}\gaa{\lceil\tfrac{n}{2} \rceil}{k+1-j}.
\end{equation*}
Another elementary computation leads to
\begin{equation*}
h(n+2,k)	-(1+q)h(n+1,k)+q\,h(n,k)-q^{2n-2k}h(n,k-1)=0.
\end{equation*}
Multiplying this by $q^{2k^2+2k}$ and summing leads to
\begin{equation*}
	\overline{T}_{n+2}-(1+q)\overline{T}_{n+1}+(q-q^{2n+2})\overline{T}_n=0,
\end{equation*}
as desired.

Additionally, we find $h(n,k)=\ga{n}{2k+1}$.

	\section{Future work}
	
	One referee has suggested to look more into the combinatorics of Durfee squares and/or generating functions/functional equations in the  context of Bressoud identities. While this  certainly looks very interesting, we feel that this would change the focus of the present paper quite a bit. 
	
	Instead, we think it would a challenge to go through Sills list~\cite{Sills} and make the present approach working for as many examples as possible. What we have done so far, is apart from the Bressoud polynomials, dealing with the Santos polynomials, 
	related with A.39/A.38 from Sills (=Slater's) list. 
	
	Here is another example, which seems to be interesting, related to A.79(=A.98). Consider the polynomials
	\begin{equation*}
		U_n=\sum_kq^{k^2}\ga{n+k}{2k}.
	\end{equation*}
	Zeilberger's algorithm produces
	\begin{equation*}
		\ga{n+2+k}{2k}-(1+q)\ga{n+1+k}{2k}-q^{2n+4-2k}\ga{n+k}{2k-2}+q\ga{n+k}{2k}=0,
	\end{equation*}
resp.
	\begin{equation}\label{recU}
		U_{n+2}-(1+q+q^{2n+3})U_{n+1}+q\,U_n=0.
	\end{equation}
	The alternative form is
	\begin{equation*}
		\overline{U}_n=\sum_jq^{15j^2-j}\ga{2n}{n-5j}-\sum_jq^{15j^2-11j+2}\ga{2n}{n+2-5j},
	\end{equation*}
	and the challenge would be to show that these polynomials also satisfy the recursion (\ref{recU}). One possible way to do that would be to split the sum into even indices $2j$ and odd indices $2j+1$. The resulting terms of the form $\ga{2n}{n+c-10j}$ could then be split with the $q$-Chu-Vandermonde into two factors of the form $\ga{n}{k+d \pm 5j}$. To work out the details might be a bit unpleasant, though.
	
\bigskip

\textbf{Acknowledgment.} Thanks are due to one referee for very substantial and constructive remarks.


\bibliographystyle{plain}

\end{document}